\documentclass[12pt]{amsart}

\title{  Multi-Hypersubstitutions and Colored Solid Varieties}
\author{K. Denecke, J. Koppitz, Sl. Shtrakov} \date{}

\def\cal{\mathcal}

\newtheorem{satz}{Theorem}[section]
\newtheorem{lemma}[satz]{Lemma}
\newtheorem{theo}[satz]{Theorem}
\newtheorem{prop}[satz]{Proposition}

\newtheorem{df}[satz]{Definition}
\newtheorem{exam}[satz]{Example}

\begin{document}


\begin{abstract}
\noindent
Hypersubstitutions are mappings which map operation symbols to terms. Terms can be visualized by trees.
Hypersubstitutions can be extended to mappings defined on sets of trees. The nodes of the trees,
describing terms, are labelled by
 operation symbols and by colors, i.e. certain positive integers. We are interested in mappings
  which map differently colored operation symbols to different terms.
In this paper we extend the theory of hypersubstitutions and solid varieties to
multi-hypersubstitutions and colored solid varieties.
We develop the interconnections between such colored terms and multi-hypersubstitutions
and the equational theory of Universal Algebra.
The collection of all varieties of a given type forms a complete lattice
which is very complex and difficult to study; multi-hypersubstitutions and colored solid varieties
   offer a new method to study complete sublattices of this lattice.

\end{abstract}

\maketitle

\section{Introduction}\label{s1}
\noindent
Let $X = \{x_1, \ldots, x_n, \ldots\}$ be a countably infinite set of
 variables, let $X_n = \{x_1, \ldots, x_n\}$ be a finite set and let $(f_i)_{i \in I}$
be a set of operation symbols where $f_i$ is $n_i-$ary. The sequence $\tau:= (n_i)_{i \in I}$
is called a type. In the usual way from variables and operation
 symbols we build up the set $W_{\tau}(X)$ of all terms of type $\tau$. An
 algebra ${\cal A} = (A; (f_i^{\cal A})_{i \in I})$ of type $\tau$ is a pair consisting
 of a set $A$ and an indexed  set of operations defined on $A$. We denote
  by $Alg(\tau)$ the class of all algebras of type $\tau$. If $s, t \in W_{\tau}(X)$, then
 the pair $s \approx t$ is called an identity in the algebra ${\cal A}$, if the term
 operations $s^{\cal A}$ and $t^{\cal A}$ induced by the terms $s$ and $t$ on
 the algebra ${\cal A}$ are equal. In this case we write ${\cal A}  \models s \approx t$.
 The binary relation $\models~ \subseteq Alg(\tau) \times W_{\tau}(X)^2$ gives
 rise to a Galois connection $(Id, Mod)$ between the power sets of $Alg(\tau)$
 and $W_{\tau}(X)^2$, where $Id$ and $Mod$ are defined for $K \subseteq Alg(\tau)$ and
 $\Sigma \subseteq W_{\tau}(X)^2$ by
\[ Id K:= \{s \approx t \mid \forall {\cal A} \in K~({\cal A} \models s \approx t)\},\] \[~ Mod
 \Sigma:= \{ {\cal A} \mid \forall s \approx t \in \Sigma~ ({\cal A} \models s \approx t)\}.\]

\noindent
As a Galois connection, $(Id, Mod)$ has the properties:

\[ \Sigma_1 \subseteq \Sigma_2 \Rightarrow Mod \Sigma_2 \subseteq Mod
\Sigma_1, ~ K_1 \subseteq K_2 \Rightarrow Id K_2 \subseteq Id K_1,\]
\[\Sigma \subseteq Id Mod \Sigma,~~ K \subseteq Mod Id K.\]
From these properties of the Galois connection $(Id, Mod)$ we obtain
that the fix points with respect to the closure operators
\[Id Mod: {\cal P}(W_{\tau}(X)^2) \to {\cal P}(W_{\tau}(X)^2) \] and
\[ Mod Id: {\cal P}(Alg(\tau)) \to {\cal P}(Alg(\tau))\]
form complete lattices
\[{\cal L}(\tau):= \{K \mid K \subseteq Alg(\tau) ~\mbox{and}~ Mod Id K = K\}\]
\[{\cal E}(\tau):= \{\Sigma \mid \Sigma \subseteq W_{\tau}(X) ^2 ~\mbox{and}~ Id Mod \Sigma = \Sigma\}\]
of all varieties of type $\tau$ and of all equational theories of type $\tau$.
These lattices are dually isomorphic.

Our next goal is to introduce two new closure operators on our sets $Alg(\tau)$ and $W_{\tau}(X)^2$
which give us complete sublattices of our two lattices ${\cal L}(\tau)$ and ${\cal E}(\tau)$.
 The new operators  are based on the concept
of hypersatisfaction of an identity by a variety. We begin with the definition of a
hypersubstitution.
A complete study of hypersubstitutions
may be found in \cite{d2}.

A hypersubstitution of type $\tau$ is a map which
associates to every operation symbol $f_i$ a term $\sigma(f_i)$ of type $\tau$,
of the same arity as $f_i$.
Any hypersubstitution $\sigma$ can be uniquely extended to a map
${\hat \sigma}$ on the set $W_{\tau}(X)$  of all terms of type $\tau$ as follows:

\noindent
\begin{enumerate}
\item[\mbox{\rm{(i)}}] If $t$ = $x_j$ for some $j \geq 1$, then $\hat\sigma[t]$ = $x_j$;
\item[\mbox{\rm{(ii)}}] if $t$ = $f_i(t_1,\ldots,t_{n_i})$ for some $n_i$-ary operation
symbol $f_i$ and some terms
$t_1$,$\ldots,t_{n_i}$, then $\hat\sigma[t]$ = $\sigma(f_i)(\hat\sigma[t_1],\ldots,
\hat\sigma[t_{n_i}])$.
\end{enumerate}

\noindent
Here the right side of (ii) means the composition of the term
$\sigma(f_i)$ and the terms $\hat\sigma[t_1], \ldots, \hat\sigma[t_{n_i}]$.\\

We can define a binary operation $\circ_h$ on the set $Hyp(\tau)$
of all hypersubstitutions of type $\tau$, by taking
$\sigma_1 \circ_h \sigma_2$ to be the hypersubstitution which
maps each fundamental operation symbol $f_i$ to the term
$\hat\sigma_1[\sigma_2(f_i)]$. That is, \\
\[ \sigma_1 \circ_h \sigma_2 := \hat\sigma_1 \circ \sigma_2, \]
where $\circ$ denotes the ordinary composition of functions. The operation $\circ_h$
is associative. The {\it identity hypersubstitution}
$\sigma_{id}$ which maps every $f_i$ to $f_i(x_1,\ldots,x_{n_{i}})$ is an identity
element for this operation.
Then ${\cal H}yp(\tau):=(Hyp(\tau); \circ_h, \sigma_{id})$ is a monoid.

\begin{df} \rm \label{d1.1}Let ${\cal M}$ be any submonoid of ${\cal H}yp(\tau)$.
An algebra ${\cal A}$ is said to
$M$-hypersatisfy
\index{hypersatisfaction}
an identity $u \approx v$ if for every hypersubstitution
$\sigma \in M$, the identity
${\hat \sigma}[u]$ $\approx$ ${\hat \sigma}[v]$ holds in ${\cal A}$. In this case we say
that the identity $u \approx v$ is an {\it $M$-hyperidentity} of ${\cal A}$ and we write
${\cal A} \begin{tabular}[t]{c}$\models$ \\[-0.9ex] {\scriptsize
\it {M-hyp}} \end{tabular}$ $u \approx v$.
For $M = Hyp(\tau)$ we write ${\cal A} \begin{tabular}[t]{c}$\models$ \\[-0.9ex] {\scriptsize
\it {hyp}} \end{tabular}$ $u \approx v$.

An identity is called an $M$-hyperidentity of a variety $V$ if it holds as an
$M$-hyperidentity in every algebra in $V$.
A variety $V$ is called {\it $M$-solid} if every identity of
$V$ is an $M$-hyperidentity of $V$.
When $M$ is the whole monoid $Hyp(\tau)$, an $M$-hyperidentity
is called a {\it hyperidentity}, and an $M$-solid variety is called
a {\it solid} variety.
\end{df}

\noindent Let ${\cal M}$ be any submonoid of ${\cal H}yp(\tau)$.
Since $M$ contains the identity hypersubstitution, any
$M$-hyperidentity of a variety $V$ is an identity of $V$. This means
that the relation of $M$-hypersatisfaction, defined between
$Alg(\tau)$ and $W_{\tau}(X)^2$, is a subrelation of the relation of
satisfaction from which we induced our Galois-connection $(Id,
Mod)$. The new Galois-connection induced by the relation of
$M$-hypersatisfaction is called  $(H_MMod, H_MId)$ and is defined
 on classes $K$ and sets
$\Sigma$ as follows:
\\
\centerline{$ H_MId K =$} \\ \centerline{$   \{ s \approx t \in
W_{\tau}(X)^2 \mid s \approx t ~\mbox{is~
an}~M-~\mbox{hyperidentity~ of~} {\cal A}~ \mbox{for ~all}~ {\cal A}
\in K \},$}\\
\centerline{$H_MMod \Sigma =$} \\ \centerline{$  \{ {\cal A} \in Alg(\tau)\mid ~\mbox{all~ identities~ in~} \Sigma~\mbox{are}~
~ M-~\mbox{hyperidentities~ of~} {\cal A}  \}.$}\\

\noindent
The Galois-closed classes of algebras under this connection are the $M$-solid
 varieties of type $\tau$, which form a complete sublattice of the lattice of
  all varieties of type $\tau$.  Thus studying $M$-solid and
solid varieties is a way to study complete sublattices of the lattice
of all varieties of a given type. \\

We now introduce some closure operators on the two sets $Alg(\tau)$ and $W_{\tau}(X)^2$.
For equations we define
\[\chi^E_M[u \approx v]:= \{\hat\sigma[u] \approx \hat\sigma[v] \mid \sigma \in M\}.\] \\
\vspace*{-0.2cm}

\noindent
For any set $\Sigma$ of identities we set

\[\chi^E_M[\Sigma] = \bigcup \limits_{u \approx v \in \Sigma}\chi^E_M[u \approx v].\]

\noindent
Hypersubstitutions can also be applied to algebras, as follows. Given an
algebra ${\cal A}$ = $(A; (f_i^{\cal A})_{i \in I})$
and a hypersubstitution $\sigma$, we define the algebra $\sigma({\cal A})$ =
 $(A; (f_i^{\sigma({\cal A})})_{i \in I})$: = $(A; (\sigma(f_i)^{\cal A})_{i \in I})$.
This algebra is called the {\it derived algebra}
determined by ${\cal A}$ and $\sigma$. Notice that by definition
it is of the same type as the algebra ${\cal A}$. Now we define an operator
$\chi^A_M$ on the set $Alg(\tau)$,
first on individual algebras and then on classes $K$ of algebras, by
\[ \chi^A_M[{\cal A}] = \{\sigma[{\cal A}] \mid \sigma \in M \},\] and
\[\chi^A_M[K] = \bigcup \limits_{{\cal A} \in K}\chi^A_M[{\cal A}].\]
\noindent
If $M = Hyp(\tau)$ the operators are denoted by $\chi^A$ and $\chi^E$.

Let $\tau$ be a fixed type and let ${\cal M}$ be any submonoid of ${\cal H}yp(\tau)$.
The two operators $\chi^E_M$ and $\chi^A_M$ are  closure operators and
are connected by the condition
\[ \chi^A_M[{\cal A}] \mbox{ satisfies } u \approx v \mbox{  iff  }
{\cal A} \mbox{ satisfies } \chi^E_M [u \approx v]. \]

\noindent
The following propositions are also obvious (see \cite{d1} or  \cite{d2}).

\begin{theo}\label{t1.2}
 Let $K\subseteq Alg(\tau )$ and $\Sigma \subseteq W_\tau(X)$.
Then there holds\\

\begin{tabular}{llll}
\rm{(i)}&$H_MMod\Sigma$& = &$Mod\chi_M ^{E}[\Sigma],$ \\
\rm{(ii)}& $H_MMod\Sigma$& $\subseteq$& $Mod\Sigma,$\\
\rm{(iii)}& $\chi_M^{A}[H_MMod\Sigma ]$&=&$H_MMod\Sigma,$\\
\rm{(iv)}& $\chi_M ^{E}[IdH_MMod\Sigma ]$&=&$IdH_MMod\Sigma, $\\
\rm{(v)}& $H_MModH_MIdK$&=&$ModId\chi_M ^{A}[K],$\\

\rm{(i')}&$H_MIdK$&=&$Id\chi_M ^{A}[K],$\\
\rm{(ii')}&$H_MIdK$&$\subseteq$&$ IdK,$\\
\rm{(iii')}& $\chi_M ^{E}[H_MIdK]$&=&$H_MIdK,$\\
\rm{(iv')}& $\chi_M^{A}[ModH_MIdK]$&=&$ModH_MIdK,$\\
\rm{(v')}& $H_MIdH_MMod\Sigma$& =&$IdMod\chi_M ^{E}[\Sigma].$\\

\end{tabular}
\end{theo}

\noindent
$M$-solid varieties can be characterized by the following theorem:

\begin{theo}\label{t1.3}
 Let ${\cal M}$ be a monoid of hypersubstitutions of type $\tau$. For any variety $V$ of type
$\tau$, the following conditions are equivalent: \\

\begin{tabular}{llll}
\rm{(i)}&$V$ &=& $H_MMod H_MId V,$\\
\rm{(ii)}&$\chi^A_M[V]$& =& $V,$\\
\rm{(iii)}&$Id V$ &=& $H_MId V,$\\
\rm{(iv)}& $\chi^E_M[Id V]$ &=& $Id V$. \\
\end{tabular}
\vspace*{1cm}

\noindent
\rm{And dually, for any equational theory $\Sigma$ of type $\tau$, the
following conditions are equivalent:}\\

\begin{tabular}{llll}
\rm{(i')}&$\Sigma$&=& $H_MId H_MMod \Sigma,$\\
\rm{(ii')}&$\chi^E_M[\Sigma]$ &=& $\Sigma,$\\
\rm{(iii')}&$Mod \Sigma$& =& $H_MMod \Sigma,$\\
\rm{(iv')}& $\chi^A_M[Mod \Sigma]$ &=& $Mod \Sigma$.
\end{tabular}
\end{theo}
\hspace*{4cm}~~~~~~~~~~~~~~\hfill $\rule{2mm }{2mm}$

\noindent
A variety which satisfies the condition (i) is called an $M$-
hyperequational class.

\noindent
The subrelation
$\begin{tabular}[t]{c}$\models$ \\[-0.9ex] {\scriptsize
\it {M-hyp}} \end{tabular}$ of $\models$ satisfies the following condition:\\
$\forall K \subseteq Alg(\tau)~ \forall \Sigma \subseteq W_{\tau}(X)^2~((H_MMod
 \Sigma = K~ \wedge~ H_MId K = \Sigma) \Rightarrow Mod \Sigma = K ~\wedge~ Id K = \Sigma).$

Such subrelations are called Galois-closed. The Galois connections induced by
 Galois-closed subrelations of a given relation generate complete sublattices of
 the complete lattices generated by the Galois connection induced by that given relation  \cite{d3}.
Moreover we have the following result.

\begin{theo}\label{t1.4}
Let ${\cal M}$ be a monoid of hypersubstitutions of type $\tau$. Then the class ${\cal S}_M(\tau)$
of all $M$-solid varieties of type $\tau$ forms a complete sublattice of the lattice ${\cal L}(\tau)$
of all varieties of type $\tau$.
Dually, the class of all $M$-hyperequational theories forms a complete sublattice of the lattice
of all equational theories of type $\tau$. \hfill \qed
\end{theo}

\noindent
When ${\cal M}_1$ and ${\cal M}_2$ are both submonoids of ${\cal H}yp(\tau)$
and ${\cal M}_1$ is a submonoid of
${\cal M}_2$, then the corresponding complete lattices satisfy ${\cal S}_{M_2}(\tau)$ $\subseteq$
${\cal S}_{M_1}(\tau)$. As a special case, for any ${\cal M}$ $\subseteq {\cal H}yp(\tau)$
 we see that the lattice
${\cal S}(\tau)$ of all solid varieties of type $\tau$ is always a sublattice of the
lattice ${\cal S}_M(\tau)$.
At the other extreme, for the smallest possible submonoid ${\cal M}$ = $\{\sigma_{id}\}$
the corresponding
lattice of $M$-solid varieties is the whole lattice ${\cal L}(\tau)$ of all varieties of
 type $\tau$. Thus
we obtain a range of complete sublattices of ${\cal L}(\tau)$ to ${\cal S}({\tau})$.

Our aim is to transfer this theory to another kind of hypersubstitution and to colored terms.

\section{Multi-hypersubstitutions and colorations of terms}\label{s2}

In a term a certain operation symbol may occur more than once. To distinguish
between the different occurrences of the same operation symbol we assign to each
 occurrence of any operation symbol a color.
Representing a term by a tree, we get a vertex-colored graph. To apply different
 hypersubstitutions to the same operation symbol, if it is differently colored,
 we define the concept of a multi-hypersubstitution.\\

\begin{df} \rm \label{d2.1}
A map $\rho $ from $\mathbb N$ into $Hyp(\tau )$ is called a multi-
hypersubstitution. Let $Hyp^{\mathbb N}$ be the set of all
multi-hypersubstitutions.
\end{df}

\noindent
To distinguish between different occurrences of the same operation symbol in a
 term $t$ we assign to each operation symbol in $t$ an address in $t$, i.e. an
  element of a given set. Usually adresses are sequences of natural numbers. Let
   $add(t)$ be the set of all addresses of the term $t$.
We introduce the concept of a coloration to allow that equal operation symbols
 at different places are considered differently. On the other hand, a coloration
  allows that equal operation symbols with different addresses are equally colored.
   A coloration of a type $\tau$ is defined in the following way.

\begin{df} \rm \label{d2.2}
Any mapping $\alpha_t$ from $add(t)$, $t \in W_{\tau}(X)\setminus
X$, into $\mathbb N$ is called a coloration of the term $t$. We denote
by $C(t)$ the set of all colorations of the term $t$. A set $C
\subseteq \bigcup \{C(r) \mid r \in W_{\tau}(X)\}$ with $|C(r) \cap
C | = 1 $
for all $r \in W_{\tau}(X)\setminus X$ is called a
coloration of $W_{\tau}(X)$.
\end{df}

\noindent Using  colorations, multi- hypersubstitutions can be
extended to mappings defined on terms. In a first step we extend
multi- hypersubstitutions to mappings from the set $Sub(t)$ of all
subterms of a given term $t$ to the set of terms.

\begin{df} \rm \label{d2.3}
Let $C$ be a coloration of $W_{\tau}(X)$, $t \in W_{\tau}(X)$ with the
coloration $\alpha \in C$, $s \in Sub(t)$, and let $\rho$ be a multi-hypersubstitution.
\begin{enumerate}
\item[\mbox{\rm{(i)}}] If $s\in X$ then $\widehat{\rho}_{C,t}[s]:=s$.
\item[\mbox{\rm{(ii)}}] If $s=f_{i}(s_{1},\ldots,s_{n_{i}})$ with $i\in I$ and
$s_{1},\ldots,s_{n_{i}}\in Sub(t)$, where $f_i$ has the address $a$ in $t$
 then \[\widehat{\rho}_{C,t}[s]:=\rho(\alpha(a))(f_i)
 (\widehat{\rho}_{C,t}[s_{1}],\ldots,\widehat{\rho}_{C,t}[s_{n_{i}}]),\].
\end{enumerate}
\end{df}

\noindent
Using the mappings $\widehat{\rho}_{C,t}[t] : Sub(t) \rightarrow W_{\tau}(X)$
 for terms $t \in W_{\tau}(X)$ we can extend multi-hypersubstitutions to mappings defined on terms.

\begin{df} \rm \label{d2.4}
Let $C$ be a coloration of $W_{\tau}(X)$ and $\alpha \in C$ be the coloration
 of a term $t \in W_{\tau}(X)$. Then for a multi-hypersubstitution $\rho$ we
  put $\widehat{\rho}_C[t]:=\widehat{\rho}_{C,t}[t]$.
\end{df}

\noindent
The following example shows that the composition of two multi-
hypersubstitutions does not be a multi- hypersubstitution. For this
we consider the type $\tau =(2)$, the terms $s=f(y,f(y,x))$ and
$t=f(f(x,y),y)$ withe the following coloration $C$:\\
\\
$\alpha_{t}(a) = 0$ for all $a\in add(t);$\newline $\alpha_{s}(q) =
0$, where $s$ is the address of the leftmost $f$ in $s$;\newline
$\alpha_{s}(a) = 1$ for all $a\in add(s)$, $a\neq q$.\\ \\
 Let
$\rho \in Hyp(\tau)^{\mathbb N}$ be a multi-hypersubstitution with
$\rho (0)=\sigma _{yx}$ and $\rho (a)=\sigma _{id}$ for $a\in
\mathbb{N}\backslash \ \{0\}$. Then we have $\widehat{\rho}_C[t]= s$
and $\widehat{\rho}_C[t]= f(f(y,x),y)$. Since all operation symbols
in $t$ have the same color, we can not find a
 multi-hypersubstitution which provides
$f(f(y,x),y)$ by application on $t$.

If all addresses of a term have the same color, then a multi-
hypersubstitution can be replaced on that term by one of its
components, i.e. by an ordinary hypersubstitution as the following
lemma shows.

\begin{lemma}\label{l2.5}
Let $C$ be a coloration of $W_{\tau }(X)$, $\rho \in Hyp(\tau)^{\mathbb N}$, and $t\in
W_{\tau }(X)$ such that there is an $n \in \mathbb N$ with $\alpha _{t}(a)=n$
for all $a\in add(t)$. Then $\widehat{\rho }_{C}[t]=\widehat{\rho (n)}[t]$.
\end{lemma}

\noindent
{\bf Proof:}
Since $\widehat{\rho}_C[t]=\widehat{\rho}_{C,t}[t]$ we show by induction that
$\widehat{\rho (n)}[s]:=\widehat{\rho}_{C,t}[s]$ for each $s \in Sub(t)$.
If $s\in X$ then  $\widehat{\rho }_{C,t}[s]=s=\widehat{\rho (n)}[s]$.\newline
Let $s=f_{i}(s_{1},\ldots,s_{n_{i}})$ and suppose that $\widehat{\rho }
_{C,t}[s_{j}]=\widehat{\rho (n)}[s_{j}]$ for $1\leq j\leq n_{i}$ then\\
\\
$\widehat{\rho }_{C,t}[f_{i}(s_{1},\ldots,s_{n_{i}})]$

\begin{tabular}{lll}
\hspace*{1.5cm}&=&$\rho(\alpha_{t}(a))(f_{i})($ $\widehat{\rho}_{C,t}[s_{1}],\ldots,$ $\widehat{
\rho }_{C,t}[s_{n_{i}}])$\\
\hspace*{1.5cm}& & ($a$ denotes the address of $f_{i}$ in $t$)\\
\hspace*{1.5cm}&=&$\rho (\alpha _{t}(a))(f_{i})(\widehat{\rho (n)}[s_{1}],\ldots,$ $\widehat{%
\rho (n)}[s_{n_{i}}])$ \\
\hspace*{1.5cm}& & (by the hypothesis)\\
\hspace*{1.5cm}&=&$\rho (n)(f_{i})(\widehat{\rho (n)}[s_{1}],\ldots,$ $\widehat{\rho (n)}[s_{n_{i}}])$\\
\hspace*{1.5cm}&=&$\widehat{\rho (n)}[f_{i}(s_{1},\ldots,s_{n_{i}})]$.
 \hspace*{4,4cm}~~~~~~~~~~~~~~~~~~~~~~~~~~~~~~~~~\hfill $\rule{2mm}{2mm}$
\end{tabular}

\section{Colored solid Varieties}\label{s3}

\noindent
Using multi-hypersubstitutions
we define operators corresponding to $\chi_M^A, \chi_M^E$ of the
 introduction and apply these operators to sets of equations and to classes of algebras.

\begin{df} \rm \label{d3.1}
Let $C$ be a coloration of $W_{\tau}(X)$ and $\Sigma \subseteq W_{\tau}(X)^2$.
 Then we put $\chi_{C}^{e}[\Sigma]:= \{\widehat {\rho}_C[u]\approx
  \widehat{\rho}_C[v]\ |\ u\approx v\in
 \Sigma ,\ \rho \in Hyp(\tau )^{\mathbb N}\}= \chi _{C}^{e,0}[\Sigma ].$
For $n\in \mathbb N$, we put $\chi_{C}^{e,n+1}[\Sigma ]:=\chi _{C}^{e,}[\chi_{C}^{e,n}[\Sigma ] ].$
Let $\chi _{C}^{E}[\Sigma ]:= \bigcup \{\chi _{C}^{e,n}[\Sigma ]\ |\ n\in \mathbb N\}$.
\end{df}

\noindent
For any set $\Sigma \subseteq W_{\tau}(X)^2$ we have $\chi _{C}^{e}[\Sigma ]=\Sigma$
iff $\chi _{C}^{E}[\Sigma ]=\Sigma$.
Indeed, If $\chi _{C}^{e}[\Sigma ]=\Sigma$ then it is easy to see that $\chi _{C}^{e,n}[\Sigma ]
=\Sigma$ for all $n\in \mathbb
N$ and thus $\chi _{C}^{E}[\Sigma ]=\Sigma$. Conversely, if $\chi_{C}^{E}[\Sigma ]=\Sigma$
then $\chi _{C}^{e}[\Sigma ]
\subseteq \chi _{C}^{E}[\Sigma ]$ provides $\chi_{C}^{e}[\Sigma] \subseteq \Sigma$.
Since $\rho_{id} \in Hyp(\tau
)^{\mathbb N}$ we have $\Sigma \subseteq \chi _{C}^{e}[\Sigma ] $ and altogether,
 $\chi _{C}^{e}[\Sigma ]=\Sigma$.

\begin{df} \rm \label{d3.2}
Let $C$ be a coloration of $W_{\tau}(X)$ and let $V$ be a variety of type $\tau
$.
$V$ is called $C-colored\ solid$ if $IdV=\chi _{C}^{E}[IdV]$.
\end{df}

\begin{exam}\label{e3.3}
\rm{We consider the following example of a $C$-colored variety of type
 $\tau=(2)$. Let $RB = Mod\{x(yz) \approx (xy)z \approx xz, x^2 \approx x\}$
  be the variety of rectangular bands.
It is well-known that $RB$ is solid (see e.g. \cite{d6}). The set $Id RB$ of
all identities satisfied in $RB$ is the set of all equations $s \approx t$
such that the first variable of  $s$ agrees with the first variable of $t$
and the last variable of $s$ agrees with the last variable of $t$. For $t \in W_{(2)}(X)$
 such that $t$ starts and ends with the same variable we define $\alpha_t:add(t) \to
\mathbb N$ with $\alpha_t(a) = 1$ for all $a \in add(t)$ and if $t$ starts and ends with
different variables, then we define $\alpha_t(a)= 2$ for all $a \in add(t)$. Let
$s \approx t \in Id RB$ and let $\rho$ be a multi-hypersubstitution. Then it is easy to
see that $\widehat{{\rho}_C}[s] \approx \widehat{{\rho}_C}[t] \in Id RB$ using solidity.}
\end{exam}

\begin{df} \rm \label{d3.4}
Let $C$ be a coloration of $W_{\tau}(X)$,  let $\rho \in $ $Hyp(\tau )^{\mathbb N}$
and let ${\cal A}=(A;(f_{i}^{\cal A})_{i\in I})$ be an
algebra of type $\tau.$ Then we define:\newline
$\rho \lbrack {\cal A}]:=(A;{(f_{i}^{\rho[{\cal A}]})}_{i\in I})$ where
$f_{i}^{\rho[{\cal A}]} =\widehat{\rho }_C[f_{i}(x_{1},\ldots,x_{n_{i}})]^{\cal A}$ for $i\in I$.
\end{df}

\begin{df} \rm \label{d3.5}
Let $C$ be a coloration of $W_{\tau}(X)$ and $K$ be a class of algebras of
type $\tau $. Then we put $\chi _{C}^{A}[K]:=\{\rho \lbrack {\cal A}]\
|\ {\cal A}\in K,\ \rho \in $ $Hyp(\tau )^{\mathbb N}\}$.
\end{df}

\noindent
It is easy to check that $\chi _{C}^{E}$ and $\chi _{C}^{A}$ have the properties of a
completely additive closure operator. From Lemma \ref{l3.6} below will follow that
a $C$-colored solid variety is solid.
Our new closure operators are connected with the operators defined in the Introduction.

\begin{lemma}\label{l3.6}
Let $C$ be a coloration of $W_{\tau}(X)$. Then
\begin{enumerate}
\item[\mbox{\rm{(i)}}] $\chi ^{E}[\Sigma ]\subseteq \chi _{C}^{E}[\Sigma ]$ \ for each $%
\Sigma \subseteq W_\tau(X)^2$.
\item[\mbox{\rm{(ii)}}] $\chi _{C}^{A}[K]=\chi ^{A}[K]$ \ for each $K\subseteq Alg(\tau).$
\end{enumerate}
\end{lemma}
\noindent
{\bf Proof}:
(i) $\chi _{C}^{e}[\Sigma ] \subseteq \chi _{C}^{E}[\Sigma ]$ we have to show
that $\chi ^{E}[\Sigma ]\subseteq \chi
_{C}^{e}[\Sigma ]$. Let $\sigma \in Hyp(\tau )$. Then we consider the
multi-hypersubstitution $\rho \in $ $Hyp(\tau )^{\mathbb N}$ with $\rho(a)=\sigma $
 for all $a \in \mathbb N$. If \ $t\in W_{\tau }(X)$\ then we
will check that $\widehat
{\sigma}[t]=\widehat{\rho}_C[t]$. Since $\widehat{\rho}_C[t]=\widehat{\rho}_{C,t}[t]$
 we show that $\widehat
{\sigma}[s]:=\widehat{\rho}_{C,t}[s]$ for each $s \in Sub(t)$ by induction on the
complexity of the term $s$.

If $s\in X$ then $\widehat{\sigma}[s]=s=\widehat{\rho }_{C,t}[s]$.

Assume that \ $s=f_{i}(s_{1},\ldots,s_{n_{i}})$ with $i\in I$ and $s_{1},\ldots,s_{n_{i}}%
\in Sub(t)$ and  suppose inductively  that $\widehat{\sigma}[s_{k}]=\widehat{%
\rho}_{C,t}[s_{k}]$ for 1 $\leq k\leq n_{i}$. Then for
a suitable $a \in \mathbb N$ we have
\[\widehat{\rho}_{C,t}[s]=\rho (a)(f_{i})(\widehat{\rho }_{C,t}[s_{1}], \ldots,\widehat{
\rho}_{C,t}[s_{n_{i}}])\] \[=\sigma (f_{i})(\widehat{\sigma
}[s_{1}], \ldots,\widehat{
\sigma}[s_{n_{i}}])=\widehat{\sigma}[f_{i}(s_{1},\ldots,s_{n_{i}})].\]
\noindent (ii) Let ${\cal A}\in K$. For $i\in I$ let $a_i$ be the
color of $f_i$ in the fundamental term
$f_{i}(x_{1},\ldots,x_{n_{i}})$.

Let $\rho \in $ $Hyp(\tau )^{\mathbb{N}}$. Then we consider the
hypersubstitution $\sigma \in Hyp(\tau )$ with $\sigma (f_{i})=\rho(a_i)(f_{i})$ for $i\in I$ and have
\[f_{i}^{\rho \lbrack
{\cal A}]}=\widehat{\rho }_C[f_{i}(x_{1},\ldots,x_{n_{i}})]^{\cal A}=
(\rho(a_i)(f_{i})(x_{1},\ldots,x_{n_{i}}))^{\cal A}\]\[=\sigma
(f_{i})(x_{1},\ldots,x_{n_{i}})^{\cal A}=\sigma (f_{i})^{\cal A}.\] This shows that $\chi _{C}^{A}[
{\cal A}]\subseteq \chi ^{A}[{\cal A}].$

Let $\sigma \in Hyp(\tau )$. Then we consider the multi-hypersubstitution $%
\rho \in $ $Hyp(\tau )^{\mathbb N}$ with $\rho (a)=\sigma $ for all $a\in
\mathbb N$. Now we have \\
$$f_{i}^{\sigma \lbrack {\cal A}]}= \sigma
(f_{i})^{\cal A}=(\sigma (f_{i})(x_{1},\ldots,x_{n_{i}}))^{\cal A} =(\rho(a_{i})
(f_{i})(x_{1},\ldots,x_{n_{i}}))^{\cal A}$$
$$=\widehat{\rho}_C[f_{i}(x_{1},\ldots,x_{n_{i}})]^{\cal A}=f_{i}^{\rho \lbrack
{\cal A}]}.$$
 This shows that $\chi ^{A}[{\cal A}]\subseteq \chi
_{C}^{A}[{\cal A}].$
Altogether we have $\chi _{C}^{A}[{\cal A}]=\chi ^{A}[{\cal A}]$
and thus $\chi _{C}^{A}[K]=\chi ^{A}[K].$
\hfill $\rule{2mm}{2mm}$\\

\noindent
Using the operators $\chi_{C}^{A}$ and $\chi_{C}^{E}$
we define two new relations between $Alg(\tau)$ and $W_{\tau}(X)^2$.

\begin{df} \rm \label{d3.7}
Let $C$ be a coloration of $W_{\tau}(X)$. Then we put \\
\\
$R_{1}:=\{({\cal A},s\approx t)\ |\ {\cal A}\in Alg (\tau ),\
s\approx t\in W_{\tau}(X)^2 ,\ \chi _{C}^{E}[s\approx t]\subseteq Id{\cal A}%
\ \};$\\ \\
$R_{2}:=\{({\cal A},s\approx t)\ |\ {\cal A}\in Alg (\tau ),\
s\approx t\in W_{\tau}(X)^2 ,\ s\approx t\in Id\chi _{C}^{A}[{\cal A}]\ \};$%
\\ \\
$\mathcal{C}Mod\Sigma :=\{\ {\cal A}\ |\ \{\ {\cal A}\ \}\times
\Sigma \subseteq R_{1}\}$ for $\Sigma \subseteq W_{\tau}(X)^2;$\\ \\
$\mathcal{C}IdK:=\{s\approx t\ |\ K\times \{\ s\approx t\ \}\subseteq
R_{1}\} $ for $K\subseteq Alg (\tau ).$
\end{df}

\noindent
Because of Lemma  \ref{l3.6}  (ii) the relation $R_2$ agrees with the relation
$\begin{tabular}[t]{c}$\models$ \\[-0.9ex] {\scriptsize
\it {hyp}} \end{tabular}$.
Since
\begin{tabular}[t]{c}$\models$ \\[-0.9ex] {\scriptsize
\it {hyp}} \end{tabular} is a Galois closed subrelation of $\models$,
 the relation $R_2$ has the same property. Since we are more interested
 in $R_1$, we ask whether $R_1$ is a Galois closed subrelation of $\models$.
At first we have the following property:

\noindent
\begin{prop}\label{p3.8}
Let $C$ be a coloration of $W_{\tau }(X)$. Then $R_{1}$ is a subrelation of $
\models $ such that for $K\subseteq A\lg (\tau )$ and all $\Sigma \subseteq
W_\tau(X)^2 $ the following holds: If $\ \Sigma =\mathcal{C}IdK$ and $K=\mathcal{C}Mod\Sigma $
then $K=Mod\Sigma .$
\end{prop}

\noindent
{\bf Proof}:
At first we show that $R_{1}$ is a subrelation of $\models$. For this let $(%
{\cal A},s\approx t)\in R_{1}.$ Then $\chi _{C}^{E}[s\approx t]\subseteq Id {\cal A}$
 and $s\approx t\in Id{\cal A}$ since $
\chi _{C}^{E}$ is a closure operator. Thus $({\cal A},s\approx t)\in \, \models$.

Now we show that $\Sigma = \chi _{C}^{E}[\Sigma ]$. Since $\chi
_{C}^{E}$ is a closure operator, we have $\Sigma \subseteq \chi
_{C}^{E}[\Sigma ].$

Conversely let $u\approx v\in \chi _{C}^{E}[\Sigma ].$ Then there is an $%
s\approx t\in \Sigma $ with $u\approx v\in \chi _{C}^{E}[s\approx t]$, i.e.
 $\chi _{C}^{E}[u\approx v]\subseteq \chi _{C}^{E}[\chi _{C}^{E}[s\approx
t]]\subseteq \chi _{C}^{E}[s\approx t].$ From $s\approx t\in \Sigma =
\mathcal{C}IdK$ it follows $\chi _{C}^{E}[s\approx t]\subseteq IdK$ and $
\chi _{C}^{E}[u\approx v]\subseteq \chi _{C}^{E}[s\approx t]\subseteq IdK$,
 i.e. $u\approx v\in \mathcal{C}IdK.$ This shows $\chi _{C}^{E}[\Sigma
]\subseteq \Sigma $.
Now we have \\
\\ \centerline{
\begin{tabular}{lll}
$Mod\Sigma
$&=&$\{\ {\cal A}\ |\ {\cal A}\in Alg (\tau ),\Sigma \subseteq Id%
{\cal A}\ \}$\\ &&\\
&=&$\{\ {\cal A}\ |\ {\cal A}\in Alg (\tau ),\ \chi _{C}^{E}[\Sigma
]\subseteq Id{\cal A}\ \}$\\ &&\\
&=&$\mathcal{C}Mod\Sigma $\\ &&\\
&=&$K$.
\end{tabular}}

~~~ \hfill $\rule{2mm}{2mm}$\\

\noindent
To obtain a characterization of colored solid varieties we check at
first the conditions of Theorem \ref{t1.2} for M-solid varieties.

\noindent
\begin{lemma}\label{l3.9}
Let $C$ be a coloration of $W_{\tau}(X)$, let $\Sigma \subseteq W_\tau(X)^2$,
 and let $K$ be a class of algebras of type $\tau$. Then
\[\chi_C^E[\Sigma] \subseteq Id K \Leftrightarrow
 \chi_{C}^{E}[\Sigma ]\subseteq Id\chi _{C}^{A}[K].\]
\end{lemma}

\noindent
{\bf Proof}:
Because of $K \subseteq \chi_C^A[K] \Rightarrow Id \chi_C^A[K] \subseteq Id K$ we
 get $\chi_C^E[\Sigma] \subseteq Id \chi_C^A[K] \Rightarrow  \chi_C^E[\Sigma] \subseteq Id K.$

By Lemma \ref{l3.6} (i) and the definition of the closure operators $\chi _{C}^{E}$ and $\chi ^{E}$
we have $\chi _{C}^{E}[\Sigma ]=\chi _{C}^{E}[\chi _{C}^{E}[\Sigma
]]\supseteq \chi ^{E}[\chi _{C}^{E}[\Sigma ]]\supseteq \chi _{C}^{E}[\Sigma
] $, i.e. $\chi _{C}^{E}[\Sigma ]=\chi ^{E}[\chi _{C}^{E}[\Sigma ]]$.

Consequently we have

$$\chi _{C}^{E}[\Sigma]\subseteq Id K
\Rightarrow \chi ^{E}[\chi_{C}^{E}[\Sigma ]]\subseteq Id K \Rightarrow
\chi_{C}^{E}[\Sigma]\subseteq Id \chi^{A}[K], $$
since $(\chi^E, \chi^A)$ is a conjugate pair.
Then by Lemma \ref{l3.6}, $\chi_C^E[\Sigma] \subseteq Id \chi_C^A[K].$
\hfill $\rule{2mm}{2mm}$\\

\noindent
As a consequence we have
\[{\cal C}Mod \Sigma = \{{\cal A} \mid {\cal A} \in Alg(\tau)~
\mbox{and}~ \chi_C^E[\Sigma] \subseteq Id \chi_C^A[{\cal
A}]\}.\quad\quad\quad
 \hfill (*)\]

\noindent
Now we prove that the sets of the form $\mathcal{C}IdK$
and $\mathcal{C}Mod\Sigma$ are closed under the operators
 $\chi _{C}^E$ and $\chi _{C}^A$, respectively.

\begin{lemma}\label{l3.10}
Let $C$ be a coloration of $W_{\tau}(X)$.
\begin{enumerate}
\item[\mbox{\rm{(i)}}] For $K\subseteq A\lg (\tau )$ there holds $\chi _{C}^{E}[\mathcal{C}IdK]=
\mathcal{C}IdK.$
\item[\mbox{\rm{(ii)}}] For $\Sigma \subseteq W_\tau(X)^2$ there holds $\chi _{C}^{A}[\mathcal{C}
Mod\Sigma ]=\mathcal{C}Mod\Sigma .$
\end{enumerate}
\end{lemma}

\noindent
{\bf Proof:}
(i) Clearly, $\mathcal{C}IdK\subseteq \chi _{C}^{E}[\mathcal{C}IdK]$.
 Let $u\approx v\in \chi _{C}^{E}[\mathcal{C}IdK]$. Then there is
an equation $s\approx t\in \mathcal{C}IdK$ with $u\approx v\in \chi _{C}^{E}[s\approx
t]$. Since $\chi _{C}^{E}$ is a closure operator we have $\chi _{C}^{E}[u\approx v]
\subseteq \chi _{C}^{E}[\chi _{C}^{E}[s\approx t]]=\chi
_{C}^{E}[s\approx t]$. From $s\approx t\in \mathcal{C}IdK$ it follows $\chi
_{C}^{E}[s\approx t]\subseteq IdK$. Then $\chi _{C}^{E}[u\approx v] \subseteq Id K$,
thus $u \approx v \in
{\cal C}IdK$.\\
(ii) Clearly, $\mathcal{C}Mod\Sigma \subseteq \chi _{C}^{A}[\mathcal{C}Mod\Sigma ]$.
Let ${\cal A}\in \chi _{C}^{A}[\mathcal{C}Mod\Sigma ]$. Then there is
a ${\cal B}\in \mathcal{C}Mod\Sigma $ such that ${\cal A}\in \chi
_{C}^{A}[{\cal B}]$. This implies $\chi _{C}^{A}[{\cal A}]\subseteq \chi _{C}^{A}
[\chi _{C}^{A}[{\cal B}]]=\chi _{C}^{A}[{\cal B}]$ and
$Id\chi _{C}^{A}[{\cal B}]\subseteq Id\chi _{C}^{A}[{\cal A}]$.
 From ${\cal B}\in \mathcal{C}Mod\Sigma $ it follows $
\chi _{C}^{E}[\Sigma ]\subseteq Id\chi _{C}^{A}[{\cal B}]$ by $(*)$.
Thus $\chi _{C}^{E}[\Sigma ]\subseteq Id\chi _{C}^{A}[{\cal A}]$ and
${\cal A}\in \mathcal{C}Mod\Sigma $ by $(*)$.
This shows that $\chi_{C}^{A}[\mathcal{C}Mod\Sigma]\subseteq \mathcal{C} Mod\Sigma $.
\hfill$\rule{2mm}{2mm}$\\

Because of $R_1 \subseteq $$~\models$, we have
$\mathcal{C}Mod\Sigma \subseteq Mod\Sigma $ and
 $\mathcal{C} IdK\subseteq IdK$ for all
  sets $\Sigma \subseteq W_\tau(X)^2$ and for all $K \subseteq Alg(\tau)$.\\

The following example shows that $R_1$ is not a Galois closed subrelation of $\models$.
We consider the set of all identities of the greatest
solid variety $V_{HS}:= Mod\{x_1(x_2x_3) \approx (x_1x_2)x_3, x_1^2
 \approx x_1^4, x_1x_2x_1x_3x_1x_2 x_1 \approx x_1x_2x_3x_2x_1, x_1^2x_2^2x_3
  \approx x_1^2x_2x_1^2x_2x_3, x_1x_2^2x_3^2 \approx x_1x_2x_3^2x_2x_3^2\}$ of semigroups (\cite{d5}).

\begin{exam}\label{e3.11}
\rm{Let $\tau =(2)$ and let $f$ be the binary operation symbol.
 Then there is a coloration $C$ of $W_{(2)}(X)$ such that
$R_{1}$ is not  a Galois-closed subrelation of $\models$. We
will show that there is a set $\Sigma$ of equations with
\[\Sigma = CId V_{HS} ~\mbox{and}~ V_{HS} = C Mod \Sigma \not \Rightarrow \Sigma = Id V_{HS}.\]

\noindent
We put $s:=f(f(x,x),f(f(x,x),f(x,x)))$ and $\alpha_{s}: add(s)\longrightarrow
\mathbb{N}$ with $\alpha_{s}(a) = 1$ for all $a\in add(s).$\newline
We set $\Psi :=(\{s\approx t\ |\ t\in W_{(2)}(X)\ \setminus \{s\}\ \}\cup
\{t\approx s\ |\ t\in W_{(2)}(X)\ \setminus \{s\}\ \})\cap IdV_{HS}$ and $%
\Sigma :=IdV_{HS}$ $\setminus \ \Psi $.\newline
For $t\in W_{(2)}(X)\ \setminus \ \{s\}$ we define $\alpha_{t}:add(t)
\longrightarrow \mathbb{N}$ with $\alpha_{t}(a) = 0$ for all $a\in add(t).$\newline
This gives us a coloration $C$ of $W_{(2)}(X)$.
\newline
If $t\in W_{(2)}(X)\ \setminus \ \{s\}$ then for $\rho \in
Hyp(2)^{\mathbb{N}}$ holds $\widehat{\rho }[t]=\widehat{\rho (0)}[t]$ by Lemma \ref{l2.5}

\noindent
As a consequence we have (together with Lemma \ref{l3.6})
 $\chi _{C}^{e}[\Theta ] = \chi ^{E}[\Theta ]$ for each set $\Theta
\subseteq W_{(2)}(X)^2.$ Since $\chi ^{E}$ is a closure operator,
 we have the idempotence propertyand it is easy to check
that $\chi _{C}^{e,n}[\Sigma ] = \chi ^{E}[\Sigma ]$ for all $n \in \mathbb{N}$.
 This shows that $\chi _{C}^{E}[\Sigma ] = \chi
^{E}[\Sigma ]$. \newline
Hence $\mathcal{C}Mod\Sigma =\{\ {\cal A}\ |\ \{\ {\cal A}\
\}\times \Sigma \subseteq R_{1}\}=\{\ {\cal A}\ |\ {\cal A}\in
Alg(\tau ),\ \chi _{C}^{E}[\Sigma ]\subseteq Id{\cal A}\ \}=\{\
{\cal A}\ |\ {\cal A}\in Alg (\tau ),\ \chi ^{E}[\Sigma
]\subseteq Id{\cal A}\ \}=V_{HS}$ since $V_{HS}$ is defined by the
hyperidentity $f(f(x,y),z)\approx \ f(x,f(y,z))$.\newline

\noindent
Further we have ${\cal C}Id V_{HS}=\{r\approx t\ |\ V_{HS}\times
\{\ r\approx t\ \}\subseteq R_{1}\}$\newline
$=\{r\approx t\ |\ \ r\approx t\in W_\tau(X)^2 ,\ \chi _{C}^{E}[r\approx
t]\subseteq IdV_{HS}\}$\newline
$=\{r\approx t\ |\ \ r\approx t\in IdV_{HS},\ \chi _{C}^{E}[r\approx
t]\subseteq IdV_{HS}\}$\newline
$=\{r\approx t\ |\ \ r\approx t\in \Sigma ,\ \chi _{C}^{E}[r\approx
t]\subseteq IdV_{HS}\}\cup \{r\approx t\ |\ \ r\approx t\in \Psi ,\ \chi
_{C}^{E}[r\approx t]\subseteq IdV_{HS}\}$\newline
$=\Sigma \cup \{r\approx t\ |\ \ r\approx t\in \Psi ,\ \chi
_{C}^{E}[r\approx t]\subseteq IdV_{HS}\}$ (since $V_{HS}$ is solid ) \newline
$=\Sigma \cup \emptyset $ \newline
since the multi-hypersubstitution $\rho \in $ $%
Hyp(2)^{\mathbb{N}}$ with $\rho (0)=\sigma _{xy}$ and $\rho (a)=\sigma _{x}$
for $a\in \mathbb{N}\backslash \ \{0\}$ provides
$\widehat{\rho_C }[s]=x$
and $\widehat{\rho_C }[t]=t$ for $t\in W_{(2)}(X)\ \backslash \ \{s\},$ where $%
s\approx t\in IdV_{HS}$ implies $t\neq x$, i.e. $\widehat{\rho_C }[s]\approx
\widehat{\rho_C }[t]\notin IdV_{HS}$.\newline
\newline
Finally, because of $f(x,x)\approx s\in IdV_{HS}\ \backslash \ \Sigma $ we
have $IdV_{HS}\neq \Sigma.$
}\end{exam}

\begin{prop}\label{p3.12}
The pair $(\chi _{C}^{E}$,$\chi _{C}^{A})$ is in general not
  a conjugate pair of completely additive closure operators.
\end{prop}
\noindent
{\bf Proof}:
Assume that $(\chi _{C}^{E}$,$\chi _{C}^{A})$ forms a conjugate pair of
completely additive closure operators. Then for all ${\cal A}\in Alg
(\tau )$ and all $s\approx t\in W_\tau(X)^2$ there holds $\ s\approx t\in Id\chi
_{C}^{A}[{\cal A}]$ iff $\chi _{C}^{E}[s\approx t]\subseteq Id
{\cal A}$. In particular,  $R_{1}=\;R_{2}$. But $
R_{2}$ is a Galois-closed subrelation of $\models $  and we showed in Example \ref{e3.11}
that $R_{1}$ is not Galois-closed. Thus $R_{1}\;$and $
R_{2}$ are different, a contradiction to $R_{1}=\;R_{2}$.
\hfill $\rule{2mm}{2mm}$\\

\noindent
Since the proof of the four equivalent characterizations
 of $M$-solid varieties uses this property we cannot expect to have the same situation.

\begin{prop}\label{p3.13}
Let $C$ be a coloration of  $W_\tau(X)$, $K\subseteq Alg (\tau )$
 and $\Sigma \subseteq W_\tau(X )^2$. Then there holds
\begin{enumerate}
\item[\mbox{\rm{(i)}}] $\mathcal{C}Mod\Sigma =Mod\chi _{C}^{E}[\Sigma ]$,
\item[\mbox{\rm{(ii)}}] $\mathcal{C}IdK\subseteq Id\chi _{C}^{A}[K]$, but the converse
inclusion is in general not true.
\end{enumerate}
\end{prop}

\noindent {\bf Proof}: (i) There holds $${\cal A}\in
\mathcal{C}Mod\Sigma \quad \iff \quad\chi _{C}^{E}[\Sigma ]\subseteq
Id{\cal A}\quad \iff \quad{\cal A}\in Mod\chi _{C}^{E}[\Sigma ].$$

 \noindent (ii) Let
$\ u\approx v\in \mathcal{C}IdK$. Then $\chi _{C}^{E}[u\approx
v]\subseteq IdK$ and $\chi _{C}^{E}[u\approx v]\subseteq Id\chi
_{C}^{A}[K]$ by Lemma \ref{l3.9}. Now we have $u\approx v\in \chi
_{C}^{E}[u\approx v]\subseteq Id\chi _{C}^{A}[K]$. Altogether this
shows that $\mathcal{C}IdK\subseteq Id\chi _{C}^{A}[K]$.\newline In
order to show that the converse direction is in general not true we
consider the type $(2)$ and the coloration $C$ of $W_{(2)}(X)$ from
Example \ref{e3.11}. Moreover let $\Sigma \subseteq W_\tau (X)^2$ be
as given in Example \ref{e3.11}. Then we have $\mathcal{C}IdV_{HS}=
\Sigma \neq IdV_{HS}=Id\chi^{A}[V_{HS}]=Id\chi_{C}^{A}[V_{HS}]$,
since we have $ V_{HS}=\chi^{A}[V_{HS}]$ ($V_{HS}$ is solid) and
$\chi ^{A}[V_{HS}] = \chi _{C}^{A}[V_{HS}]$ (Lemma \ref{l3.6}).
\hfill $\rule{2mm}{2mm}$\\

\begin{prop}\label{p3.14}
Let $C$ be a coloration of  $W_\tau(X)$, $K\subseteq Alg (\tau )$ and $
\Sigma \subseteq W_\tau(X)^2$. Then there holds:
\begin{enumerate}
\item[\mbox{\rm{(i)}}]$ \chi _{C}^{E}[Id\mathcal{C}Mod\Sigma ]\supseteq Id\mathcal{C}
Mod\Sigma $, but the converse inclusion is in general not true,
\item[\mbox{\rm{(ii)}}] $\chi _{C}^{A}[Mod\mathcal{C}IdK]=Mod\mathcal{C}IdK$.
\end{enumerate}
\end{prop}
\noindent
{\bf Proof}:
(i) The inclusion is clear, since $\chi _{C}^{E}$ is a closure
operator.\newline
In order to show that the converse direction is in general not true we
consider the type $(2)$ and the coloration $C$ of $W_{(2)}(X)$ from
Example \ref{e3.11}. Moreover let $\Sigma \subseteq W_\tau (X)^2$ and  $s\in $  $
W_{(2)}(X)$ be as given in Example \ref{e3.11} and $\rho \in $
$Hyp(\tau )^{\mathbb{N}}$ be defined by $\rho(0) = \sigma_{xy}$
and $\rho(a) = \sigma_x$ for $a \in \mathbb N \setminus\{0\}$.
Then we
get $s\approx f(x,x)\in IdV_{HS}$ and thus
$\widehat{\rho }_C[s]\approx
\widehat{\rho}_C[f(x,x)]\in \chi _{C}^{E}[IdV_{HS}]$, i.e. $x\approx x^2
\in \chi _{C}^{E}[IdV_{HS}]$. Because of $x\approx x^2
\notin IdV_{HS}$ we obtain $\chi _{C}^{E}[IdV_{HS}]\neq IdV_{HS}$, i.e. $\chi
_{C}^{E}[Id\mathcal{C}Mod\Sigma ]\neq Id\mathcal{C}Mod\Sigma $ since $V_{HS} = {\cal C} Mod \Sigma$.
\newline \noindent
(ii) Clearly, $Mod\mathcal{C}IdK\subseteq \chi _{C}^{A}[Mod\mathcal{C}
IdK]$.\\
For the converse inclusion let ${\cal A}\in \chi _{C}^{A}[Mod\mathcal{C}
IdK]$. Then there is a ${\cal B}\in Mod\mathcal{C}IdK$ with $
{\cal A}\in \chi _{C}^{A}[{\cal B}]$.
We want to show that $\mathcal{C}IdK\subseteq Id{\cal A}$. For this let
$u\approx v\in \mathcal{C}IdK$, i.e. $\chi _{C}^{E}[u\approx v]\subseteq IdK$.
Since $\chi _{C}^{E}$ is a closure operator we have $\chi _{C}^{E}[\chi
_{C}^{E}[u\approx v]]\subseteq IdK$ and thus $\chi _{C}^{E}[u\approx
v]\subseteq \mathcal{C}IdK$. Since ${\cal B}\in Mod\mathcal{C}IdK$ \ we
have then $\chi _{C}^{E}[u\approx v]\subseteq Id{\cal B}$.
Lemma \ref{l3.9} shows that then $\chi _{C}^{E}[u\approx v]\subseteq Id\chi _{C}^{A}[
{\cal{ B}}]$. Moreover ${\cal A}\in \chi _{C}^{A}[{\cal B}]$
implies $Id\chi _{C}^{A}[{\cal B}]\subseteq Id{\cal A}$.
Altogether we have $u\approx v\in \chi _{C}^{E}[u\approx v]\subseteq Id\chi
_{C}^{A}[{\cal B}]\subseteq Id{\cal A}$. Consequently, $\mathcal{C}
IdK\subseteq Id{\cal A}$.\newline
Finally, $\mathcal{C}IdK\subseteq Id{\cal A}$ means ${\cal A}\in
Mod\mathcal{C}IdK$. Altogether this shows the converse inclusion $\chi
_{C}^{A}[Mod\mathcal{C}IdK]\subseteq Mod\mathcal{C}IdK$.
\hfill $\rule{2mm}{2mm}$ \\

\begin{prop}\label{p3.15}
Let $C$ be a coloration of $W_{\tau }(X)$, $K\subseteq A\lg (\tau )$ and
$\Sigma \subseteq W_\tau(X)^2$. Then there holds:
\begin{enumerate}
\item[\mbox{\rm{(i)}}] $ModId\chi _{C}^{A}[K]\subseteq \mathcal{C}Mod\mathcal{C}IdK$, but
the converse inclusion is in general not true;
\item[\mbox{\rm{(ii)}}] $\mathcal{C}Id\mathcal{C}Mod\Sigma \subseteq IdMod\chi
_{C}^{E}[\Sigma ]$, but the converse inclusion is in general not true.
\end{enumerate}
\end{prop}
\noindent
{\bf Proof:}
(i) We have\\

\noindent
\begin{tabular}{llll}
$ModId\chi _{C}^{A}[K]$
&$\subseteq$&$ Mod\mathcal{C}IdK$ & by
Proposition \ref{p3.13} (ii)\\
&=&$Mod\chi _{C}^{E}[\mathcal{C}IdK]$ & by
Lemma \ref{l3.10} (i)\\
&=&$\mathcal{C}Mod\mathcal{C}IdK$ & by Proposition \ref{p3.13} (i).\\
\end{tabular}

\noindent
In order to show that the converse direction is in general not true
we consider the type $\tau = (2)$. Further let $\beta :W_{\tau }(X)\longrightarrow
\mathbb{N}$ be a bijection. (Such a bijection exists since $W_{\tau }(X)$
is countable.) For $t\in W_{(2)}(X)$ we define $\alpha_{t}:add(t)\longrightarrow
\mathbb{N}$ with $\alpha_{t}(a)=\beta (t)$ for all $a\in add(t).$ Then
$C:=\{ \alpha_{t}\ |\ t\in W_{(2)}(X)\ \}$ is a coloration of $W_{(2)}(X)$.\\
Let $u\approx v\in W_\tau(X)^2$. If $u=v$ then obviously
$\chi _{C}^{E}[u\approx v]\subseteq \{w\approx w\ |\ w\in W_{(2)}(X)\}\subseteq IdV_{HS}$.
If $u\neq v$ and $u,v\in X$ then we get $\chi
_{C}^{E}[u\approx v]=\{u\approx v\}\not\subseteq  IdV_{HS}$. If $u\neq v$
and $u\notin X$ then $\beta (u)\neq \beta (v)$ and we may consider the
 multi-hypersubstitution $\rho \in $ $Hyp(2)^{\mathbb{N}}$ with
  $\rho (\beta (u))=\sigma _{xy}$ and $\rho (\beta (v))=\sigma _{x}$. Then $
\widehat{\rho}_C[u]=u$ and $\widehat{\rho}_C[v]=r$, where $r$ is the first
letter in $v$. Thus $\widehat{\rho}_C[u]\approx \widehat{\rho}_C[v]\notin
IdV_{HS}$, i.e. $\chi _{C}^{E}[u\approx v]\not\subseteq IdV_{HS}$.
If $u\neq v$ and $v\notin X$ then we get $\chi _{C}^{E}[u\approx
v]\not\subseteq IdV_{HS}$ in the dual manner.
Altogether we have $\mathcal{C}IdV_{HS}=\{w\approx w\ |\ w\in W_{(2)}(X)\}$
and consequently,  $\mathcal{C}Mod\mathcal{C}IdV_{HS} =
 Mod\{w\approx w\ |\ w\in W_{(2)}(X)\}=A\lg (2)$ by Lemma \ref{l3.10}
 (i) and Proposition \ref{p3.13} (i).\\
On the other hand we have $$ModId\chi _{C}^{A}[V_{HS}] = ModId\chi
^{A}[V_{HS}] = ModIdV_{HS} = V_{HS}$$ by Lemma 3.6 and since
$V_{HS}$ is solid. This shows that $$ModId\chi _{C}^{A}[V_{HS}]\neq
\mathcal{C}Mod\mathcal{C}IdV_{HS}.$$

\noindent (ii) We have $$\mathcal{C}Id\mathcal{C}Mod\Sigma \subseteq
Id\chi _{C}^{A}[\mathcal{C}Mod\Sigma ] =
 Id\mathcal{C}Mod\Sigma = IdMod\chi _{C}^{E}[\Sigma]$$
 by Proposition \ref{p3.13} (ii), Lemma \ref{l3.10} (ii)
  and Proposition \ref{p3.13} (i).  \\
\noindent
In order to show that the converse direction is in general not true we
consider the type $(2)$ and the coloration $C$ of $W_{(2)}(X)$ from
Example \ref{e3.11}. Moreover let $\Sigma \subseteq W_\tau(X )^2$ be defined as in
Example \ref{e3.11}.  Then we have $\mathcal{C}Id\mathcal{C}Mod\Sigma
=\mathcal{C} IdV_{HS}=\Sigma $.\\
On the other hand we have $IdMod\chi _{C}^{E}[\Sigma ] =
Id\mathcal{C}Mod\Sigma = IdV_{HS}\neq \Sigma $ by Proposition
\ref{p3.13} (i). Consequently, $\mathcal{C}Id\mathcal{C}Mod\Sigma \neq IdMod\chi_{C}^{E}[\Sigma ]$.
$\hfill \rule{2mm}{2mm}$\\

\section{$C$-colored hyperequational theories}\label{s4}

\begin{df} \rm \label{d4.1}
Let $\Sigma$ be an equational theory of type $\tau$.
Then $\Sigma$ is said to be a $C$-colored hyperequational
theory if $\mathcal{C}Mod\Sigma = Mod \Sigma$.
\end{df}

\noindent
$C$-colored hyperequational theories can be characterized in
 the same way as usual hyperequational theories:

\begin{theo}\label{t4.2}
Let $\Sigma  \in {\cal E}(\tau)$. Then the following
statements are equivalent:
\begin{enumerate}
\item[\mbox{\rm{(i)}}] $\mathcal{C}Mod\Sigma =Mod\Sigma $.
\item[\mbox{\rm{(ii)}}] $\Sigma =\mathcal{C}Id\mathcal{C}
Mod\Sigma $.
\item[\mbox{\rm{(iii)}}] $\Sigma =\chi _{C}^{E}$ $[\Sigma ]$.
\item[\mbox{\rm{(iv)}}]  $\mathcal{C}IdMod\Sigma = \Sigma $.
\end{enumerate}
\end{theo}
\noindent {\bf Proof:} (i) $\Longrightarrow$ (ii): Since
$\mathcal{C}Id\mathcal{C}Mod$ is a closure operator we have $\Sigma
\subseteq \mathcal{C}Id\mathcal{C}Mod\Sigma$. Conversely,
$\mathcal{C}Id\mathcal{C}Mod\Sigma \subseteq Id\mathcal{C}Mod\Sigma
=IdMod\Sigma = \Sigma$.

(ii) $\Longrightarrow$ (iii): Clearly, $\Sigma \subseteq \chi _{C}^{E}$
$[\Sigma]$.
Conversely, let $u\approx v\in \chi _{C}^{E}$ $[\Sigma ]=\chi _{C}^{E}$ $[%
\mathcal{C}Id\mathcal{C}Mod\Sigma]$ (because of (ii)).
Then there is an $s\approx t\in \mathcal{C}Id\mathcal{C}Mod\Sigma $
 with $u\approx v\in \chi
_{C}^{E}$ $[\ s\approx t]$. Because of $s\approx t\in \mathcal{C}Id\mathcal{C}Mod\Sigma $
 we have $\chi_{C}^{E}$ $[\ s\approx t]\subseteq Id\mathcal{C}
Mod\Sigma $. With $$K:= \mathcal{C}Mod \Sigma\ \mbox{  and   }\
 \mathcal{C}Id K =
 \mathcal{C}Id\mathcal{C}Mod\Sigma = \Sigma$$ by Proposition \ref{p3.8}
 we get $K = Mod\Sigma$.   Thus $Id\mathcal{C}Mod\Sigma =IdMod\Sigma =\Sigma $
  since $\Sigma \in {\cal E}(\tau)$.
Altogether we have $$u\approx v\in  \chi _{C}^{E} [s\approx
t]\subseteq Id\mathcal{C}Mod\Sigma =\Sigma.$$ This shows that $\chi
_{C}^{E}[\Sigma]\subseteq \Sigma $.

 (iii) $\Longrightarrow$ (i):
$\mathcal{C} Mod\Sigma \subseteq Mod\Sigma$ is clear. Conversely,
let ${\cal A}\in Mod\Sigma $. Then $\Sigma \subseteq Id{\cal A}$.
 Because of (iii) we have then $\chi _{C}^{E}$ $[\Sigma
]\subseteq Id{\cal A}$, i.e. ${\cal A}\in \mathcal{C}Mod\Sigma $.
This shows that $Mod\Sigma \subseteq \mathcal{C}Mod\Sigma.$

 (i)
$\Leftrightarrow$ (iv):\  Suppose that $\mathcal{C}IdMod\Sigma
=\Sigma $. Then  we put
$$K:=\mathcal{C}Mod\mathcal{C}IdMod\Sigma.$$ Now we have
$$\mathcal{C}Mod\Sigma =\mathcal{C}Mod\mathcal{C}IdMod\Sigma = K$$ and
$$\mathcal{C}IdK=\mathcal{C}Id\mathcal{C}Mod\mathcal{C}IdMod\Sigma =
\mathcal{C}IdMod\Sigma $$ (since $(\mathcal{C}Id,\mathcal{C}Mod)$ is
a Galois-connection). By Proposition \ref{p3.8} we get $K=Mod\Sigma
$, i.e.
 $\mathcal{C}Mod\mathcal{C}IdMod\Sigma =Mod\Sigma $, $\mathcal{C}Mod\Sigma =Mod\Sigma $
and $\Sigma $ is a $C$-colored hyperequational theory.

 Suppose that
$\Sigma $ is a $C$-colored hyperequational theory.  We use that
$\Sigma = IdMod\Sigma \supseteq \mathcal{C}IdMod\Sigma.$ Now let
$u\approx v\in \Sigma .$ Assume that $u\approx v\notin
\mathcal{C}IdMod\Sigma $. Then $\chi _{C}^{E}[$ $u\approx
v]\not\subseteq IdMod\Sigma $ and in particular there is an ${\cal
A}\in Mod\Sigma $ with
 $\chi_{C}^{E}[u\approx v]\not\subseteq Id{\cal A}$. Since $\Sigma $ is a $
C $-colored hyperequational theory we have $Mod\Sigma = \mathcal{C}Mod\Sigma $,
 i.e. ${\cal A}\in \mathcal{C}Mod\Sigma $ and $\chi
_{C}^{E}[\Sigma ]\subseteq Id{\cal A}$. Because of $u\approx v\in
\Sigma $ we have $\chi _{C}^{E}[u\approx v]\subseteq \chi
_{C}^{E}[\Sigma ]$ and thus $\chi _{C}^{E}[u\approx v]\subseteq
Id{\cal A}$, a contradiction. Consequently, $u\approx v\in
\mathcal{C}IdMod\Sigma $. This shows that $\Sigma \subseteq
\mathcal{C}IdMod\Sigma $. $\hfill \rule{2mm}{2mm}$

 When $\Sigma$ is
a $C$- colored hyperequational theory, then $Mod\Sigma =\chi
_{C}^{A}[Mod\Sigma]$. To\ \  see this, suppose that $Mod\Sigma \neq
\chi _{C}^{A}[Mod\Sigma ]$,\  then $Mod\Sigma \neq \chi
^{A}[Mod\Sigma ]$ by Lemma \ref{l3.6}, i.e. $Mod\Sigma $ is not
solid. Then there is an ${\cal A}\in Mod\Sigma $ with $\chi ^{E}$
$[\Sigma]\not\subseteq Id{\cal A}$.
 Because of $\chi ^{E}[\Sigma ]\subseteq
\chi _{C}^{E}$ $[\Sigma ]$ we have $\chi _{C}^{E}$ $[\Sigma ]\not\subseteq Id{\cal A}$, i.e.
${\cal A}\notin \mathcal{C}Mod\Sigma $ contradicts
$\mathcal{C}Mod\Sigma =Mod\Sigma $.\\

The following example shows that the converse is not true, i.e.
$C$-colored hyperequational theories cannot be characterized by the
 condition $Mod\Sigma = \chi_ {C}^{E}[Mod\Sigma]$.

\begin{exam}\label{e4.3}
\rm{Let $\tau =(2).$ We use the coloration $C$ from Example \ref{e3.11}.
 Further we take $\Sigma :=IdV_{HS}$. Then we have $ModIdV_{HS}=V_{HS}=\chi ^{A}[V_{HS}]$,
 since $V_{HS}$ is solid.
By Lemma \ref{l3.6} we have
$V_{HS} = \chi _{C}^{A}$ $[V_{HS}] = \chi _{C}^{A}$ $[ModIdV_{HS}]$.\\

Now we show that $\mathcal{C}ModIdV_{HS}\neq V_{HS}=ModIdV_{HS}$.
Actually we will show that $\mathcal{C}ModIdV_{HS}\subseteq {\bf B}$,
where ${\bf B}$ denotes the variety of bands. We consider the term $s$
 and the hypersubstitution $\rho $ from  Example \ref{e3.11}  and note
 that the application of $\rho $ to the      identity $s\approx f(x,x)\in $
 $IdV_{HS}$ provides the idempotent law. Thus $x\approx x{{}^2}
\in \chi _{C}^{E}$ $[IdV_{HS}]$. Let ${\cal A}\in \mathcal{C}%
ModIdV_{HS} $. Then $\chi _{C}^{E}$ $[IdV_{HS}]\subseteq Id{\cal A}$,
and in particular $x\approx x {{}^2} \in Id{\cal A}$, i.e. ${\cal A}\in {\bf B}$.}
\end{exam}

\section{Characterizations of colored solid varieties}\label{s5}

In Section \ref{s3} colored solid varieties $V$ were defined by the
property $Id V = \chi_C^E[Id V]$. We get the following characterization:

\begin{theo}\label{t5.1}
Let $C$ be a coloration of $W_{\tau }(X)$ and $V$ be a variety
of type $\tau $. Then the following statements (i)-(iii) are equivalent:
\begin{enumerate}
\item[\mbox{\rm{(i)}}] $\mathcal{C}ModIdV=V.$
\item[\mbox{\rm{(ii)}}] $IdV=\mathcal{C}IdV.$
\item[\mbox{\rm{(iii)}}] $\chi _{C}^{E}[IdV]=IdV$.
\end{enumerate}
Further, each of the statements (i)-(iii) implies both $V=\mathcal{C}Mod
\mathcal{C}IdV$ and $\chi _{C}^{A}[V]=V$, but the converse implications are
in general not true.
\end{theo}
{\bf Proof:} (i)$ \Rightarrow $(ii) Suppose that
$\mathcal{C}ModIdV=V$. Then $IdV$ is a $C$-colored hyperequational
theory, i.e. Mod$IdV=\mathcal{C}ModIdV$ and
$\mathcal{C}Id\mathcal{C}ModIdV$ $=IdV$. Thus $\mathcal{C}IdV$
$=IdV$, i.e. $V$ is $C$-colored solid.

 (ii)$ \Rightarrow $(iii)
Since $\chi _{C}^{E}$ is a closure operator we
 have $Id V \subseteq \chi _{C}^{E}$ $[Id V]$.\newline
Conversely, let $u\approx v\in \chi _{C}^{E}$ $[IdV]$. Then there is an $\
s\approx t\in IdV$ with $u\approx v\in
\chi_{C}^{E}[\ s\approx t]$. From $s\approx t\in IdV$ it follows that
 $s\approx t\in \mathcal{C}IdV$ by (ii), i.e. $
\chi _{C}^{E}$ $[\ s\approx t]\subseteq IdV$. Altogether we have
$u\approx v\in IdV$. This shows that $\chi _{C}^{E}[IdV]\subseteq
IdV$.

 (iii) $\Rightarrow $(i) $IdV$ is an equational theory. Thus we
can use Theorem \ref{t4.2}
  and from $IdV=\chi _{C}^{E}$ $[IdV]$ it follows that $\mathcal{C}ModIdV=ModIdV$
  and $IdV=\mathcal{C}Id\mathcal{C}ModIdV$. This
gives
$$IdV=\mathcal{C}Id\mathcal{C}ModIdV=\mathcal{C}IdModIdV=\mathcal{C}IdV.$$

Moreover we have $$V\subseteq
\mathcal{C}Mod\mathcal{C}IdV=\mathcal{C}ModIdV = ModIdV=V.$$ This
shows that $V=\mathcal{C}ModIdV$.

Suppose that $IdV=\mathcal{C}IdV$.  We show that
$V=\mathcal{C}Mod\mathcal{C}IdV$. First we have that $ V\subseteq
\mathcal{C}Mod\mathcal{C}IdV$ since $\mathcal{C}Mod\mathcal{C}Id$ is
a closure operator. Conversely we get that
$\mathcal{C}ModIdV\subseteq ModIdV=V$. Using $IdV=\mathcal{C}IdV$
one obtains $\mathcal{C}Mod\mathcal{C} IdV\subseteq V$.

Now we show that $V=\chi _{C}^{A}[V]$. Assume that $V\neq \chi _{C}^{A}[V]$.
Then $V \not= \chi ^{A}[V]$ by Lemma \ref{l3.6}, i.e. $V$ is not solid and
 $\chi^{E}[IdV]\not\subseteq IdV$. Because of $\chi ^{E}[IdV]\subseteq \chi_{C}^{E}[IdV]$
 we have that $\chi _{C}^{E}[IdV]\not\subseteq IdV$. This shows
that $IdV\not\subseteq \mathcal{C}IdV$, contradicting
$IdV=\mathcal{C}IdV$.

Finally we prove that the opposite implication is not satisfied. Let
$\tau =(2).$ Then there are a coloration $C$ of $W_{(2)}(X)$, a
variety $ V $ of type $(2)$ with $V=\mathcal{C}Mod\mathcal{C}IdV$
and $V=\chi _{C}^{A}[V] $ such that $\mathcal{C}IdV\neq IdV$.
Indeed, we take the coloration $C$ and the set $\Sigma $ from the
proof of Example \ref{e3.11}. There we have shown that
$\mathcal{C}Mod\Sigma =V_{HS}$ and $\mathcal{C}IdV_{HS}=\Sigma $. So
we have $V_{HS}$ $=\mathcal{C}Mod\Sigma =
\mathcal{C}Mod\mathcal{C}IdV_{HS}$.
 On the other hand there holds $V_{HS}=\chi ^{A}[V_{HS}]$ since $V_{HS}$ is
solid. Thus $V_{HS}=\chi _{C}^{A}[V_{HS}]$ by Lemma \ref{l3.6}.
Moreover there holds $\mathcal{C}IdV_{HS}=\Sigma \neq IdV_{HS}$.
\hfill $\rule{2mm}{2mm}$ \\

$C$-colored hyperequational theories also have the following properties.

\begin{prop}\label{p5.2}
Let $C$ be a coloration of $W_{\tau }(X)$ and $\Sigma $ be an
equational theory. If $Id\mathcal{C}Mod\Sigma =\Sigma $ then $\Sigma
$ is a $C$-colored hyperequational theory. The converse direction is
not true.
\end{prop}
\noindent
{\bf Proof}:
Suppose that $Id\mathcal{C}Mod\Sigma =\Sigma $. Clearly, $\Sigma \subseteq
\mathcal{C}Id\mathcal{C}Mod\Sigma $. Let $u\approx v\in \mathcal{C}Id\mathcal{C}Mod\Sigma $.
 Then $\chi _{C}^{E}[u\approx v]\subseteq Id\mathcal{C}Mod\Sigma $, i.e.
 $u\approx v\in \chi _{C}^{E}[u\approx v]\subseteq Id\mathcal{C}Mod\Sigma =\Sigma $
  and $u\approx v\in \Sigma $. This shows $\mathcal{C}Id\mathcal{C}Mod\Sigma \subseteq \Sigma $.
  Altogether we have $ \Sigma =\mathcal{C}Id\mathcal{C}Mod\Sigma $. By Theorem \ref{t4.2} we have $
\mathcal{C}Mod\Sigma =Mod\Sigma $, i.e. $\Sigma $ is a $C$-colored
hyperequational theory.

 In order to show that the converse
direction is not true we consider the type $(2)$ and the coloration
$C$ of $W_{(2)}(X)$ from Example \ref{e3.11}. Moreover let $\Sigma
\subseteq W_\tau(X)^2$ be defined as in Example \ref{e3.11}. Then we
have $\Sigma
=\mathcal{C}IdV_{HS}=\mathcal{C}Id\mathcal{C}Mod\Sigma$, i.e.
$\Sigma $ is a $C$-colored hyperequational theory. But there holds $
Id\mathcal{C}Mod\Sigma =IdV_{HS}\neq \Sigma $. \hfill
$\rule{2mm}{2mm}$

\begin{prop}\label{p5.3}
Let $C$ be a coloration of $W_{\tau }(X)$ and $K$ be a class of algebras of
type $\tau $. Then $Mod\mathcal{C}IdK=\mathcal{C}Mod\mathcal{C}IdK$.
\end{prop}
\noindent
{\bf Proof:}
We have $\mathcal{C}Mod\mathcal{C}IdK$
$=\{{\cal A}\ |\ {\cal A}\in A\lg (\tau ),\ \chi _{C}^{E}[\mathcal{C}IdK]\subseteq Id{\cal A}\ \}
=\{{\cal A}\ |\ {\cal A}\in A\lg (\tau ),\ \mathcal{C}
IdK\subseteq Id{\cal A}\ \} = Mod\mathcal{C}IdK$ by Lemma 3.10 (i).
\hfill $\rule{2mm}{2mm}$

\section{Examples}\label{s6}
\noindent
 For a variety to be $C$-colored solid, it must satisfy all
identities
 obtained by applying all multi-hypersubstitutions to all identities of
 the variety. This can be difficult to verify. In the hyperidentity case,
  if we want to check whether a variety   of the form $V = Mod \Sigma$ is
   solid, we have to apply the hypersubstitutions only to the set $\Sigma$.
   This is based on the following theorem:

\begin{theo} \cite{d1}\label{t6.1}
Let $K\subseteq A\lg (\tau )$ and $\Sigma \subseteq W_\tau(X)^2$ and let
 ${\cal M} \subseteq {\cal H}yp(\tau)$ be a submonoid. Then the following  holds:
\begin{enumerate}
\item[\mbox{\rm{(i)}}]$\chi^{E}_M[\Sigma ]\subseteq IdMod\Sigma \Longleftrightarrow
IdMod\Sigma =H_MIdH_MMod\Sigma $.
\item[\mbox{\rm{(ii)}}]$\chi^{E}_M[\Sigma ]\subseteq IdMod\Sigma
\Longleftrightarrow \chi^{E}_M[IdMod\Sigma ]\subseteq IdMod\Sigma $.
\item[\mbox{\rm{(iii)}}] $\chi^{A}_M[K]\subseteq ModIdK\Longleftrightarrow ModIdK=H_MModH_MIdK$.
\item[\mbox{\rm{(iv)}}] $\chi^{A}_M[K]\subseteq ModIdK\Longleftrightarrow \chi^{A}_M[ModIdK]\subseteq
ModIdK$.
\end{enumerate}
\end{theo}

\noindent For multi-hypersubstitutions and colored terms we have:

\begin{theo}\label{t6.2}
Let $K\subseteq A\lg (\tau )$, let $C$ be a coloration of
$W_{\tau}(X)$ and $\Sigma \subseteq W_\tau(X)^2$. Then the following
holds:
\begin{enumerate}
\item[\mbox{\rm{(i)}}]$\chi_ {C} ^{E}[\Sigma ]\subseteq IdMod\Sigma \Longleftrightarrow
IdMod\Sigma  \supseteq \mathcal{C}Id\mathcal{C}Mod\Sigma $.
\item[\mbox{\rm{(ii)}}]$\chi_ {C} ^{E}[\Sigma ]\subseteq IdMod\Sigma \Longleftarrow \chi_{C}
^{E}[IdMod\Sigma ]\subseteq IdMod\Sigma $.
\item[\mbox{\rm{(iii)}}] $\chi_{C} ^{A}[K]\subseteq ModIdK\Longleftarrow ModIdK=
\mathcal{C}Mod\mathcal{C}IdK$.
\item[\mbox{\rm{(iv)}}] $\chi_{C} ^{A}[K]\subseteq ModIdK\Longleftrightarrow \chi_{C}
^{A}[ModIdK]\subseteq ModIdK$.
\end{enumerate}
The converse implications are in general not satisfied.
\end{theo}
{\bf Proof:} (i): Suppose that $\chi _{C}^{E}[\Sigma ]\subseteq
IdMod\Sigma $. Then we have $$\mathcal{C}Id\mathcal{C}Mod\Sigma
\subseteq IdMod\chi _{C}^{E}[\Sigma ] \subseteq IdModIdMod\Sigma =
IdMod\Sigma $$ by Proposition \ref{p3.15} (i).

 Suppose that
$IdMod\Sigma \supseteq \mathcal{C}Id\mathcal{C}Mod\Sigma $. Then we
have $$\chi _{C}^{E}[\Sigma ] \subseteq
\mathcal{C}Id\mathcal{C}Mod\chi _{C}^{E}[\Sigma ] =
\mathcal{C}Id\mathcal{C}Mod\Sigma \subseteq IdMod\Sigma$$ because of
${\cal C}Mod \Sigma  = Mod \chi_C^E[\Sigma] = Mod
\chi_C^E[\chi_C^E[\Sigma]] = {\cal C}Mod\chi_C^E[\Sigma]$
(Proposition \ref{p3.13}).

 In order to show that the implication
$$\chi _{C}^{E}[\Sigma ]\subseteq IdMod\Sigma \Longrightarrow
IdMod\Sigma =\mathcal{C}Id\mathcal{C}Mod\Sigma $$ is in general not
true we consider the type $(2)$ and the coloration $C$ of $
W_{(2)}(X)$ from Example \ref{e3.11}. Moreover let $\Sigma \subseteq
W_\tau(X)^2$ be defined as in Example \ref{e3.11}. Then we have
\[\chi _{C}^{E}[\chi _{C}^{E}[\Sigma ]]=\chi _{C}^{E}[\Sigma ]\subseteq
IdMod\chi _{C}^{E}[\Sigma ].\]

On the other hand we have
$$\mathcal{C}Id\mathcal{C}Mod\chi _{C}^{E}[\Sigma
]=\mathcal{C}Id\mathcal{C} Mod\Sigma = \mathcal{C}IdV_{HS}=\Sigma $$
and
$$IdMod\chi _{C}^{E}[\Sigma ]=Id\mathcal{C}Mod\Sigma = IdV_{HS}\neq \Sigma
.$$ This shows that indeed $\chi _{C}^{E}[\chi _{C}^{E}[\Sigma
]]\subseteq IdMod\chi _{C}^{E}[\Sigma ]$ but $IdMod\chi
_{C}^{E}[\Sigma ]\neq \mathcal{C}
 Id\mathcal{C}Mod\chi _{C}^{E}[\Sigma ]$.\\
(ii): This implication is clear since $\chi_{C}^E[\Sigma] \subseteq
\chi_{C}^E[IdMod\Sigma]$.

 In order to show that the converse
direction is in general not true we consider the type $(2)$ and the
coloration $C$ of $W_{(2)}(X)^2$ from Example \ref{e3.11}. Moreover
let $\Sigma \subseteq W_\tau(X)$ be defined as in Example
\ref{e3.11}. Then we have $\chi _{C}^{E}[\chi _{C}^{E}[\Sigma
]]=\chi _{C}^{E}[\Sigma ]\subseteq IdMod\chi _{C}^{E}[\Sigma ]$.

On the other hand we have $IdMod\chi _{C}^{E}[\Sigma ]=IdV_{HS}$.
Further we have $x\approx x{{}^2} \in \chi _{C}^{E}[IdV_{HS}]$.
Since $x\approx x{{}^2}\notin IdV_{HS}$ we have $\chi
_{C}^{E}[IdV_{HS}]\not\subseteq IdV_{HS}$ and thus $$\chi
_{C}^{E}[IdMod\chi _{C}^{E}[\Sigma ]]=\chi
_{C}^{E}[IdV_{HS}]\not\subseteq IdV_{HS}=IdMod\chi _{C}^{E}[\Sigma
].$$ \noindent (iii): Suppose that
$ModIdK=\mathcal{C}Mod\mathcal{C}IdK$. Then we have $\chi
_{C}^{A}[K]\subseteq ModId\chi _{C}^{A}[K] \subseteq
\mathcal{C}Mod\mathcal{C}IdK = ModIdK$ by Proposition \ref{p3.15}
(i).

In order to show that the converse implication is in general not
true we consider the type $(2)$ and the coloration $C$ of
$W_{(2)}(X)$ from Example \ref{e3.11}. \ Then we have $\chi
_{C}^{A}[V_{HS}]=V_{HS} \subseteq ModIdV_{HS}$.

 On the other hand we
have $$\mathcal{C}Mod\mathcal{C}IdV_{HS}=Mod\mathcal{C}IdV_{HS}$$ $$
 = Mod\{w\approx w\ |\ w\in W_{(2)}(X)\}=A\lg (2)\neq
 V_{HS}=ModIdV_{HS}$$
 (see proof of Proposition \ref{p3.15} (i)).

\noindent
(iv): is clear.
\hfill $\rule{2mm}{2mm}$\\

\noindent
The following example shows once more that $(\chi_{C}^E, \chi_ {C}^A)$
does not form a conjugate pair.

\begin{prop}\label{p6.3}
If \ $\tau =(n_{i})_{i\in I}$ with $n_{k}\geq 2$ for some $k\in I$ then
there are a coloration $C$ of $W_{\tau }(X)$, terms \ $s,\ t\in W_{\tau }(X)$,
 \and a multi-hypersubstitution $\rho
\in $ $Hyp(\tau )^{\mathbb{N}}$ and an algebra ${\cal A}$ of type $\tau
$ such that \[{\cal A} \models \widehat{\rho}_C[s]\approx \widehat{\rho}_C[t]\Longleftrightarrow
\rho[{\cal A}]  \models s \approx t \] is not valid.
\end{prop}
{\bf Proof}: We will show that for the free algebra ${\cal
F}_{\tau}(X)$ of type $\tau $ over an alphabet $X$ there are a
coloration $C$ of $W_{\tau }(X) $ and some $s,t\in W_{\tau }(X)$ and
an $\rho \in $ $Hyp(\tau )^{\mathbb{N}}$ such that $\widehat{\rho
}[s]\approx \widehat{\rho }[t]\in Id{\cal F}_{\tau}(X)$ but $
s\approx t\notin Id\rho \lbrack {\cal F}_{\tau}(X)]$.

 Let $k\in I$
with $n_{k}\geq 2$. We put
$s:=f_{k}(f_{k}(x_{1},\ldots,x_{1}),x_{2},\ldots,x_{2})$ and
$t:=f_{k}(f_{k}(x_{1},\ldots,x_{1}),x_{1},\ldots,x_{1})$ and we
consider a coloration $C$ of $W_{\tau }(X)$ with the following
properties:\\

\begin{enumerate}
\item[\mbox{(i)}] $\alpha_{s}(j)=0$ for all $j\in add(s),$
\item[\mbox{(ii)}] $\alpha_{t}(j)= 0$ for all $j\in add(t)$,
\item[\mbox{(iii)}] $\alpha_{f_{k}(x_{1},\ldots,x_{n_{k}})}(j)=1$ for all $j\in
add(f_{k}(x_{1},\ldots,x_{n_{k}}))$.\\
\end{enumerate}

\noindent Finally, let $\rho \in $ $Hyp^{\mathbb{N}}$ given by $\rho
(0)=\sigma _{x_{1}}$ and $\rho (j)=\widetilde{\sigma }$ for $j\in
\mathbb{N}\backslash \{0\}$, where $\sigma _{x_{1}}$ is defined by
$\sigma _{x_{1}}:f_{j}\longrightarrow x_{1}$ for all\ $j\in I$ and $
\widetilde{\sigma }$ is defined by $\widetilde{\sigma
}:f_{j}\longrightarrow x_{n_{j}}$ for all $j\in I. $

Then we have $\widehat{\rho}_C[s]=\widehat{\rho}_C[t]=x_{1}.$ This
shows that $ \widehat{\rho}_C[s]\approx \widehat{\rho}_C[t]\in
Id{\cal F}_{\tau}(X). $

 On the other hand we have $$f_{k}^{\rho
\lbrack {\cal F}_{\tau}(X)]}=
\widehat{\rho}_C[f_{k}(x_{1},\ldots,x_{n_{k}})]^{{\cal
F}_{\tau}(X)}=x_{n_{k}}^{{\cal F}_{\tau}(X)}.$$
 Then the following
holds:\\

 \noindent $s^{\rho \lbrack
\cal F_{\tau}(X)]}(x_1,x_2)$

\begin{tabular}{lll}
&=&$f_{k}^{\rho \lbrack {\cal F}_{\tau}(X)]}(f_{k}^{\rho \lbrack
{\cal F}_{\tau}(X)]}(x_{1},\ldots,x_{1}),x_{2},\ldots,x_{2})$\\
&=&$x_{n_{k}}^{F_{\tau
}(X)}(x_{n_{k}}^{{\cal F}_{\tau}(X)}(x_{1},\ldots,x_{1}),x_{2},\ldots,x_{2})$\\
&=&$x_{2}$

\end{tabular}

\noindent
and\\

\noindent $t^{\rho \lbrack {\cal F}_{\tau}(X)]}(x_{1},x_{2})$

\begin{tabular}{lll}
&=&$f_{k}^{\rho
\lbrack {\cal F}_{\tau}(X)]}(f_{k}^{\rho \lbrack \cal F_{\tau
}(X)]}(x_{1},\ldots,x_{1}),x_{1},\ldots,x_{1})$\\
&=&$x_{n_{k}}^{\cal F_{\tau
}(X)}(x_{n_{k}}^{{\cal F}_{\tau}(X)}(x_{1},\ldots,x_{1}),x_{1},\ldots,x_{1})$\\
&=&$x_{1}$.\\
\end{tabular}

\noindent
This shows $s^{\rho \lbrack {\cal F}_{\tau}(X)]}\neq t^{\rho \lbrack \cal F_{\tau
}(X)]} $ and $s\approx t\notin Id\rho \lbrack {\cal F}_{\tau}(X)].$
\hfill $\rule{2mm}{2mm}$

\begin{exam}\label{e6.4}
\rm{In the next example we will determine a collection of \ colorations $C$ of $W_{(2)}(X)$ such
that all solid varieties of bands are $C$-colored solid.
 There are exactly three nontrivial solid varieties of bands: the variety $RB$ of all rectangular
bands, the variety $NB$ of all normal bands, and the variety $RegB$ of all regular bands (see \cite{d4}).

We split the set $W_{(2)}(X)$ into two sets $A$ and $B$:\\

\begin{center}
\begin{tabular}{l}
$A:=\bigcup \{W_{(2)}(\{w\})\mid w\in X\}$ and \\
$B:=W_{(2)}(X)\setminus A$.\\ \\
\end{tabular}
\end{center}

\noindent
 By $\mathcal{F}$ we denote the set of all mappings

\[\beta :\bigcup \{\{(t,a)\mid a\in add(t)\}\mid t\in A\}\rightarrow \mathbb{N}.\]
\noindent
For $\beta \in \mathcal{F}$ we define a coloration $C_{\beta }=\{\alpha
_{t}^{\beta }\mid t\in W_{(2)}(X)\}$ as follows: For $t\in W_{(2)}(X)$ and $%
a\in add(t)$ we put\\

\begin{center}
\begin{tabular}{ll} $\alpha _{t}^{\beta }(a)=1$ &
if $t\in B$; \\ $\alpha _{t}^{\beta }(a)=\beta (t,a)$&
if $t\in A$.\\ \\
\end{tabular}
\end{center}

\noindent
Then the varieties $\ RB$, $NB$, and $RegB$ are $C_{\beta }$%
-colored solid. Indeed, let $V$ be one of the varieties $RB$, $NB$,
and $RegB$ and let $s\approx t\in IdV$. Further let $\rho \in Hyp(2)^{\mathbb N}$.%

If $s,t\in A$ then there is a variable $w\in X$ such that $s,t\in W_{(2)}(\{w\})$.
(Otherwise the identity $s\approx t$ provides the identity $x\approx y$
because of the associative and the idempotent law.) Thus $\widehat{\rho }%
_{C_{\beta }}[s]$, $\widehat{\rho }_{C_{\beta }}[t]\in W_{(2)}(\{w\})$.
Because of the associative and the idempotent law $\widehat{\rho }_{C_{\beta
}}[s]\approx $ $\widehat{\rho }_{C_{\beta }}[t]$ is an identity in $V$.

If $s,t\in B$ then $\widehat{\rho }_{C_{\beta }}[s]=\widehat{\rho (1)}[s]$
and $\widehat{\rho }_{C_{\beta }}[t]=\widehat{\rho (1)}[t]$ by Lemma \ref{l2.5}.
Since $\ V$ is solid, we have $\widehat{\rho (1)}[s]\approx \widehat{\rho (1)%
}[t]\in IdV$ and thus $\widehat{\rho }_{C_{\beta }}[s]\approx \widehat{\rho }%
_{C_{\beta }}[t]\in IdV$.

If $s\in A$ and $t\in B$ then there is a $w\in X$ such that $s\in
W_{(2)}(\{w\})$. Since $\ V$ is a variety of bands we have $\widehat{\rho }%
_{C_{\beta }}[s]\approx w\approx \widehat{\rho (1)}[s]$. By Lemma \ref{l2.5} we have $%
\widehat{\rho }_{C_{\beta }}[t]=\widehat{\rho (1)}[t]$, Since $V$ is solid
we have $\widehat{\rho (1)}[s]\approx \widehat{\rho (1)}[t]\in IdV$ and thus
$\widehat{\rho }_{C_{\beta }}[s]\approx \widehat{\rho }_{C_{\beta }}[t]\in
IdV$.

If $s\in B$ and $t\in A$ then we get in a similar way that \ $\widehat{\rho }%
_{C_{\beta }}[s]\approx \widehat{\rho }_{C_{\beta }}[t]\in IdV$.}
\end{exam}

\begin{exam}\label{e6.5}
\rm{We consider the following coloration $C=\{\alpha_{t}\mid t\in W_{(2)}(X)\}$
of $W_{(2)}(X)$:  for $t\in W_{(2)}(X)$ and $a\in add(t)$ we put\\

\begin{center}\begin{tabular}{ll}
$\alpha _{t}(a)=1$ & if $t=f(x,x)$;\\
$\alpha _{t}(a)=2$ & if $t\neq f(x,x)$. \\
\end{tabular}\\ \end{center}

\noindent
The varieties $RB$, $NB$, and $RegB$ are the only nontrivial
$C$-colored solid varieties of semigroups. Indeed, by
Example \ref{e6.4}, $RB$, $NB$, and $RegB$ are $C$-colored solid.

Conversely, let $V$ be a $C$-colored solid variety of semigroups and $\rho
\in Hyp(2)^{\mathbb N}$ with
\[\rho (1)=\sigma _{x}~\mbox{and}~\rho (i)=\sigma _{xy}~\mbox{for}~2\leq i \in
\mathbb{N}.\]
Since $V$ is $C$-colored solid, $V$ is solid and thus $V\subseteq V_{HS}$,
i.e. $f(x,x)\approx f(f(x,x),f(x,x))\in IdV$. Then we have:\\

\begin{center}\begin{tabular}{lll}
\noindent
$x$&$\approx$&$ \widehat{\rho (1)}[f(x,x)]$\newline
$\approx $ $\widehat{\rho }_{C}[f(x,x)]$ (by Lemma \ref{l2.5} )\\
&$\approx $& $\widehat{\rho }_{C}[f(f(x,x),f(x,x))] (V$ is $C$-colored solid)\\
&$\approx $& $\widehat{\rho (2)}[f(f(x,x),f(x,x))]$   (by Lemma \ref{l2.5})\\
&$\approx$&$ f(f(x,x),f(x,x))$\\
&$\approx$&$ f(x,x)$.\\
\end{tabular}\\ \end{center}

\noindent
This shows that $V$ is  a variety of bands. But $RB$, $NB$,
and $RegB$ are the only nontrivial solid varieties of bands.}
\end{exam}

\noindent
Authors' Addresses:

K. Denecke\\
Universit\"at Potsdam\\
Fachbereich Mathematik\\
Postfach 601553\\
D-14415 Potsdam\\
email:{\ttfamily kdenecke@rz.uni-potsdam.de}\\

J. Koppitz\\
Universit\"at Potsdam\\
Fachbereich Mathematik\\
Postfach 601553\\
D-14415 Potsdam\\
email:{\ttfamily koppitz@rz.uni-potsdam.de}\\

Sl. Shtrakov\\
South-West-University Blagoevgrad\\
Faculty of Mathematics and Natural Sciences\\
2700 Blagoevgrad,
Bulgaria\\
e-mail:{\ttfamily shtrakov@swu.bg}

\end{document}